\documentclass[12pt]{article}
\usepackage{amsfonts,amssymb}
\makeatletter 
\def\EquationsBySection{\def\theequation
{\thesection.\arabic{equation}}%
\@addtoreset{equation}{section}}

\newcommand\old[1]{}
\newtheorem{theorem}{Theorem}[section]
\newtheorem{definition}[theorem]{Definition}

\newtheorem{remark}[theorem]{Remark}

\EquationsBySection \makeatother

\begin{document}
\pagestyle{plain}
\title
{\bf On the non-equivalence of Lorenz System and Chen
System\thanks{This research is supported by the National Natural
Science Foundation of China (10671212,90820302).}}

\author{Zhenting Hou$^1$,  Ning Kang$^{1,2}$, Xiangxing Kong$^1$, Guanrong Chen$^{3}$\\
and Guojun Yan$^{4}$\thanks{Email:
yanguojun2002@sina.com}\\
$^{1}$ School of Mathematics, Central South University,\\
Changsha, Hunan, 410075, China\\
$^{2}$ School of Mathematics and Computational Science, \\
Fuyang Teachers College, Fuyang, Anhui, 236041, China\\
$^{3}$ Department of Electronic Engineering, City University of Hong Kong,\\
Hong Kong, China\\
$^{4}$ Department of Mathematics, Zhengzhou University,\\
Zhengzhou, Henan, 450000, China }

\date{}
\maketitle

{\bf Abstract:} In this paper, we prove that the Chen system with a
set of chaotic parameters is not smoothly equivalent to the Lorenz
system with any parameters.


%\newpage

\section{\large\bf \quad Introduction}
\quad\quad Nonlinear science had experienced an unprecedented and
vigorous development particularly during the second half of the 20th
century, and it was considered ``the third revolution'' in natural
science in the history. The main subjects in the study of nonlinear
science include chaos, bifurcation, fractals, solitons and
complexity. Because of the important significance to unveil the
essence of chaos and wide potential application prospects of chaos
theory in many fields, research on chaos always carries a heavy
weight in nonlinear science. H. Poincar\'{e} [1] and C. Maxwell [2]
both had some vague concepts of chaos in their times. In the earlier
1960s, E. N. Lorenz [3] discovered the now-famous Lorenz system,
which actually produces visible chaos. Lorenz system is the first
mathematical and physical model of chaos, thereby becoming the
starting point and foundation stone for later research on chaos
theory. Since the 1960s, particularly with this model,
mathematicians, physicists and engineers from various fields have
thoroughly studied the essence of chaos, characteristics of chaotic
systems, bifurcations, routes to chaos, and many other related
topics [4]. There are also some chaotic systems of great
significance that are closely related to the Lorenz system, where a
particular example in point is the Chen system. Since the Chen
system was first found in 1999 [5,6], hundreds of papers have been
published on this new chaotic system with deep and comprehensive
results obtained. A monograph on the Lorenz systems family including
the Chen system has also been published [7]. To further understand
the interesting Chen system, one fundamental question has to be
answered: are the Chen system and the Lorenz system non-equivalent,
either topologically or smoothly? The purpose of this paper is to
prove that the Lorenz system and the Chen system are indeed
non-equivalent smoothly.

\bigbreak

\section{\large\bf \quad Results and Proofs}
\quad\quad The dynamical system $\phi^{abc}_t$ defined by
\begin{eqnarray}
   \left\{\begin{array}{ll}
    & \displaystyle\frac{dx}{dt}=a(y-x),\\
    & \displaystyle\frac{dy}{dt}=cx-xz-y,\\
    & \displaystyle\frac{dz}{dt}=xy-bz,
    \end{array}\right.
\end{eqnarray}
is called the {\bf Lorenz system} with parameters $a,b,c$.

The dynamical system  $\psi^{abc}_t$ defined by
\begin{eqnarray}
   \left\{\begin{array}{ll}
    & \displaystyle\frac{dx}{dt}=a(y-x),\\
    & \displaystyle\frac{dy}{dt}=(c-a)x-xz+cy,\\
    & \displaystyle\frac{dz}{dt}=xy-bz,
    \end{array}\right.
\end{eqnarray}
is called the {\bf Chen system} with parameters $a,b,c$.

It is clear that system (2.1) has 3 equilibrium points if
$b(c-1)>0$, i.e.,
 \begin{eqnarray*}
P_1&=&(0,0,0), \\
P_2&=& (-\sqrt{b(c-1)},-\sqrt{b(c-1)},c-1) ,\\
P_3&=& (\sqrt{b(c-1)},\sqrt{b(c-1)},c-1) ,
 \end{eqnarray*}
and system (2.2) has 3 equilibrium points if $b(2c-a)>0$, i.e.,
\begin{eqnarray*}
Q_1&=&(0,0,0), \\
Q_2&=& (-\sqrt{b(2c-a)},-\sqrt{b(2c-a)},2c-a) ,\\
Q_3&=& (\sqrt{b(2c-a)},\sqrt{b(2c-a)},2c-a) .
 \end{eqnarray*}

Denote the coordinates of $P_i$ by $(x_i,y_i,z_i)$, $i=1,2,3$,  and
the coordinates of $Q_i$ by $(x'_i,y'_i,z'_i)$,$i=1,2,3$, and denote
the vector fields on the right sides of  (2.1) and (2.2) by
$\vec{U}(x,y,z)$ and $\vec{V}(x,y,z)$, respectively. It is clear
that their Jacobians are:
 \begin{eqnarray*}
 D\vec{U}(x,y,z)&=& \left(
 \begin{array}{ccc}
 -a&a&0\\
 c-z&-1&-x\\
 y&x&-b
 \end{array}\right), \\
 D\vec{V}(x,y,z)&=& \left(
 \begin{array}{ccc}
 -a&a&0\\
 c-a-z&c&-x\\
 y&x&-b
 \end{array}\right), \\
 \end{eqnarray*}
and their determinants are:
\begin{eqnarray*}
\det D\vec{U}(P_1)&=& ab(c-1),  \\
 \det D\vec{U}(P_2) &=& \det D\vec{U}(P_3)=-2ab(c-1),
 \end{eqnarray*}

 \begin{eqnarray*}
 \det D\vec{V}(Q_1)&=&ab(2c-a),  \\
 \det D\vec{V}(Q_2)&=&\det D\vec{V}(Q_3)=-2ab(2c-a).
\end{eqnarray*}

In general, let $f(x)$ and $g(x)$ be vector fields on
$\mathbb{R}$$^n$, and
  \begin{equation}
    \dot{x}=f(x), \ x\in {\mathbb{R}}^n
 \end{equation}
 \begin{equation}
    \dot{y}=g(y), \ y\in {\mathbb{R}}^n
 \end{equation}
be two systems of differential equations on $\mathbb{R}$$^n$.

\begin{definition}
If there exists a diffeomorphism $h$ on $\mathbb{R}$$^n$ such that
\begin{equation}
 f(x)=M^{-1}(x)g(h(x)),
\end{equation}
where $M(x)$ is the Jacobian of $h$ at the point $x$, then (2.3) and
(2.4) are said to be {\bf smoothly equivalent}.
\end{definition}

\begin{remark}
If (2.3) and (2.4) are { smoothly equivalent}, and suppose that
$x_0$ and $y_0=h(x_0)$ are the corresponding equilibria of $f(x)$
and $g(x)$, $A(x_0)$ and $B(y_0)$ are the Jacobians of $f(x)$ and
$g(x)$, respectively, then $A(x_0)$ and $B(y_0)$ are similar, i.e.,
their
 characteristic polynomials and eigenvalues  are the same.
\end{remark}

\begin{theorem}
The Chen system and the Lorenz system are not smoothly equivalent,
i.e., there exists a Chen system $\psi^{a'b'c'}_t$ which is not
smoothly equivalent to any Lorenz system $\phi^{abc}_t$.
\end{theorem}

{\bf Proof} Since $M(x)\neq 0, M^{-1}(x)\neq 0$ and due to (2.5), we
have
$$
 f(x)=0\Leftrightarrow g(h(x))=0,
$$
that is, the equilibria $x$ of $f(\cdot)$ correspond to the
equilibria $h(x)$ of $g(\cdot)$, therefore a Chen system is
smoothly equivalent to a Lorenz system with the same number of
equilibria. It suffices to prove that a Chen system with 3
equilibria can not be smoothly equivalent to a Lorenz system with
any 3 equilibria.

Suppose that $h$ is a diffeomorphism on $\mathbb{R}$$^3$ such that
$\psi^{a'b'c'}_t$ and $\phi^{abc}_t$ smoothly equivalent under $h$.
Because $det D\vec{U}(P_1)$
 and $det D\vec{V}(Q_1)$ are positive, $P_1$ corresponds to $Q_1$ under $h$.
 Because
 \begin{eqnarray*}
D\vec{U}(P_1)&=& \left(
 \begin{array}{ccc}
 -a&a&0\\
 c&-1&0\\
0&0&-b
 \end{array}\right),
  \end{eqnarray*}
   \begin{eqnarray*}
 D\vec{V}(Q_1)= \left(
 \begin{array}{ccc}
 -a'&a'&0\\
 c'-a'&c'&0\\
 0&0&-b'
 \end{array}\right),
 \end{eqnarray*}
the characteristic equation of  $D\vec{U}(P_1)$ is:
$$\lambda^3+(a+b+1)\lambda^2+(a+ab-ac+b)\lambda-ab(c-1)=0,$$
and the characteristic equation of  $D\vec{V}(Q_1)$ is:
$$\lambda^3+(a'+b'-c')\lambda^2+(a'^2+a'b'-2a'c'-b'c')\lambda-a'b'(2c'-a')=0.$$
Let
\begin{eqnarray}
u &=& a'+b'-c', \\
  v &=& a'^2+a'b'-2a'c'-b'c',\\
  w&=&-a'b'(2c'-a').
\end{eqnarray}
By Remark 2.1, we must have
\begin{eqnarray}
a+b+1&=& u, \\
a+ab-ac+b&=& v,\\
-ab(c-1)&=&w.
\end{eqnarray}
By (2.9), we have $a=u-1-b$, so that in combining with (2.11),
$$c=1-\frac{w}{ab}=1-\frac{w}{b(u-1-b)}.$$
Substituting $a, c$ in (2.10), we have
\begin{eqnarray}b^3-ub^2+va-w=0.
\end{eqnarray}

It is clear that $D\vec{U}(P_2)$ and $D\vec{U}(P_3)$ have the same
characteristic equations, and $D\vec{V}(Q_2)$ and $D\vec{V}(Q_3)$
have the same characteristic equations. Hence, we may assume that
$P_2$ corresponds to $Q_2$. By comparing the coefficients of their
first-order terms, we have
$$ab+bc=b'c'.$$
Substituting them into the formulas of $a,c$, we get
\begin{eqnarray*}
&&b(u-1-b)+b\left(1-\frac{w}{b(u-1-b)}\right)\\
&=&(u-1)b-b^2+b-\frac{w}{u-1-b}\\
 &=&b'c',
\end{eqnarray*}
and
$$b^3-(2u-1)b^2+[(u-1)^2+u-1+b'c']b-w-(u-1)b'c'=0.$$
Subtracting this from (2.12), we obtain
\begin{eqnarray}
(u-1)b^2+(u+v-u^2-b'c')b+(u-1)b'c'=0.
\end{eqnarray}
By resultant elimination [8], a necessary and sufficient condition
for (2.12) and (2.13) to have same roots is

\begin{eqnarray}
&&M_0(a',b'c')\nonumber\\
 &&=\left|\begin{array}{lllll}
    1&     -u&         v&        -w&  0\\
    0&      1&        -u&         v&  -w\\
 u-1& u+v-u^2-b'c'& (u-1)b'c'&         0&  0\\
    0&    u-1&   u+v-u^2-b'c'& (u-1)b'c'&  0\\
    0&      0&       u-1&   u+v-u^2-b'c'&  (u-1)b'c'\\
    \end{array}\right|\nonumber\\
 &&=0.
\end{eqnarray}
Substituting $u,v,w$ in
 the above equation by (2.6), (2.7) and (2.8), we get an algebraic equation of
 $a',b',c'$, as
\begin{eqnarray}
&&b'(a'-2c')^2(1+c')(a'^3-a^4+a'^5+a'^2b'-a'^3b'+a'b'^2-2a'^2b'^2-2a'b'^3+a'^2b'^3+a'b'^4\nonumber\\
&&-2a'^2c'+3a'^3c'-4a'^4c'-3a'b'c'+4a'^2b'c'-5a'^3b'c'-b'^2c'+5a'b'^2c'-2a'^2b'^2c'+2b'^3c'\nonumber\\
&&-2a'b'^3c'-b'^4c'-2a'^2c'^2+4a'^3c'^2-3a'b'c'^2+6a'^2b'c'^2-b'^2c'^2+4a'b'^2c'^2+b'^3c'^2)\nonumber\\
 &&=0,
\end{eqnarray}
and its solution is given by the union of the following four
surfaces:
\begin{eqnarray*}
b'=0,\\
a'-2c'=0,\\
1+c'=0,
\end{eqnarray*}
\begin{eqnarray*}
&&a'^3-a^4+a'^5+a'^2b'-a'^3b'+a'b'^2-2a'^2b'^2-2a'b'^3+a'^2b'^3+a'b'^4-2a'^2c'+3a'^3c'\nonumber\\
&&-4a'^4c'-3a'b'c'+4a'^2b'c'-5a'^3b'c'-b'^2c'+5a'b'^2c'-2a'^2b'^2c'+2b'^3c'-2a'b'^3c'\nonumber\\
&&-b'^4c'-2a'^2c'^2+4a'^3c'^2-3a'b'c'^2+6a'^2b'c'^2-b'^2c'^2+4a'b'^2c'^2+b'^3c'^2=0.
\end{eqnarray*}
Denote the  point set of all solutions of (2.14) by $C$. It is clear
that $C$ is a Borel subset of $\mathbb{R}$$^3$ and its Lebegue
measure is 0. So, there are many points not belonging to $C$, for
example, the values of $a'=45,b'=5,c'=28$ give
$b'(2c'-a')=55>0,M_0(45,5,28)=2.919\times10^{11}\neq 0$, i.e.,
$(45,5,28)\notin C$. This means that the Chen system
$\psi_t^{45,5,28}$ is not smoothly equivalent to the Lorenz system
$\phi^{abc}_t$ with any values of $a,b,c$, while $\psi_t^{45,5,28}$
is chaotic according to [5, 6 or 7 (p.39)].

\centerline{\bf References} \bigbreak

\def\REF#1{\par\hangindent\parindent\indent\llap{#1\enspace}\ignorespaces}

\footnotesize

\REF{[1]} J. H. Poincar\'{e}, Sur le probl\'{e}me des trois corps et
les \'{e}quations de la dynamique. Divergence des s\'{e}ries de M.
Lindstedt, Acta Mathematica, vol. 13, pp. 1-270, 1890

\REF{[2]} J. C. Maxwell, Does the progress of physical sciences tend
to give any advantage to the opinion of necessity (or determinism)
over that of contingency or events and the freedom of will? In: L.
Campbell and W. Garnett (Eds.), The life of James Clerk Maxwell,
with a selection from his correspondence and occasional writings and
a sketch of his contributions to science, Macmillan, London, pp.
434-444, 1882

\REF{[3]} E. N. Lorenz, Deterministic nonperiodic flow, J. Athmosph.
Sc. vol. 20, pp. 130-141, 1963

\REF{[4]} C. Sparrow, The Lorenz Equations: Bifurcations, Chaos, and
Strange Attractors, Springer, New York, 1982

\REF{[5]} G. Chen and T. Ueta, Yet another chaotic attractor, Int.
J. of Bifurcation and Chaos, vol. 9, pp. 1465-1466, 1999

\REF{[6]} T. Ueta and G. Chen, Bifurcation analysis of Chen's
attractor, Int. J. of Bifurcation and Chaos, vol. 10, pp. 1917-1931,
2000

\REF{[7]} G. Chen and J. L\"{u}, Dynamical Analysis, Control and
Synchronization of the Lorenz Systems Family, Science Press,
Beijing, 2003

\REF{[8]} L. Yang, J. Zhang and X. Hou, Nonlinear Systems of
Algebraic Equations and Mechanical Theorem Proving, Shanghai
Scientific and Technological Education Publishing House, Shanghai,
1996

\end{document}